\documentclass[12pt]{amsart}
\usepackage{latexsym, amssymb, amsmath}

\newtheorem{theorem}{Theorem}[section]
\newtheorem{lemma}[theorem]{Lemma}

\newtheorem{remark}[theorem]{Remark}
\theoremstyle{definition}

\newtheorem{example}{Example}

\newtheorem{claim}{Claim}

\makeatletter
    
    \@addtoreset{equation}{section}
  \makeatother

\newcommand{\Ric}{{\rm Ric}}

\newcommand{\n}{\nabla}

\begin{document}

\title[Rotational Symmetry of Shrinking Gradient Yamabe Solitons]
{Rotational symmetry of complete shrinking gradient Yamabe solitons}

\author{Shun Maeta}
\address{Department of Mathematics, Faculty of Education, and
Department of Mathematics and Informatics, Graduate School of Science and Engineering, 
Chiba University, 1-33, Yayoicho, Inage, Chiba, 263-8522, Japan.}
\curraddr{}
\email{shun.maeta@gmail.com~{\em or}~shun.maeta@chiba-u.jp}
\thanks{The author is partially supported by the Grant-in-Aid for Scientific Research (C), No.23K03107, Japan Society for the Promotion of Science.}


\date{}

\dedicatory{}

\keywords{shrinking gradient Yamabe solitons; Yamabe soliton version of the Perelman conjecture; rotationally symmetric}

\commby{}

\begin{abstract}
In this paper, we show that any nontrivial complete shrinking gradient Yamabe soliton whose scalar curvature is bounded below by the soliton constant everywhere and is strictly greater than the constant at some point is rotationally symmetric.
This assumption is optimal for higher dimensions.
This result resolves the Yamabe-soliton analogue of Perelman's conjecture. 

\end{abstract}

\maketitle

\bibliographystyle{amsplain}


\section{History and the result}\label{intro}

Geometric flows are fundamental and highly effective tools for analyzing the structure of Riemannian manifolds. 
In particular, the Ricci flow and the Yamabe flow, both introduced by R. Hamilton, are widely regarded as prototypical geometric flows (see \cite{Hamilton82}, \cite{Hamilton89}).
Over the past few decades, both the Ricci flow and the Yamabe flow have developed rapidly and have become central topics in geometry (see \cite{Brendle05}, \cite{Brendle07}, \cite{BS08}, \cite{Chow92}, \cite{CLN06}, \cite{Hamilton95}, \cite{SS03}, \cite{Ye94}).

To understand geometric flows, it is important to study their singularity models. Ricci solitons 
$
\Ric + \nabla \nabla F=\lambda g,
$
and Yamabe solitons 
$
\nabla \nabla F=(R- \lambda) g,
$
are special solutions of these flows, and are expected to serve as singularity models (see \cite{CLN06}).
Here, $R$ and $\Ric$ denote the scalar curvature and the Ricci tensor of a Riemannian manifold $(M,g)$, respectively, $\n\n F$ is the Hessian of a smooth function $F$ on $M$, and $\lambda$ is a constant.
From this perspective, these solitons have been studied extensively over the last few decades.

As is well known, Perelman conjectured that any three-dimensional complete noncompact $\kappa$-noncollapsed gradient steady Ricci soliton with positive curvature is rotationally symmetric (see \cite{Perelman02}). 
This is known as Perelman's conjecture.
After the conjecture was proposed, significant partial results were obtained (see \cite{Brendle11}, \cite{CC12}, \cite{CCCMM14}, \cite{CMM16}, \cite{WW13}, and the references therein). 
Finally, S. Brendle \cite{Brendle13} completely resolved the conjecture.
However, in higher dimensions, the problem is still not well understood, and considerable research efforts are currently devoted to its investigation (see, for example, \cite{Brendle14},  \cite{CMM17}).

It is also important to understand the structure of gradient shrinking Ricci solitons (see \cite{Chowetal07}, \cite{Perelman02}). 
In low dimensions, the structure of gradient shrinking Ricci solitons is well known (see \cite{BM15}, \cite{Hamilton88}, \cite{Hamilton89}, \cite{Perelman02}, \cite{CCZ08}, \cite{NW08}, \cite{Wu91}). 
For higher dimensions, under certain curvature assumptions (in particular, on the Riemannian curvature tensor), gradient shrinking Ricci solitons have been studied (see \cite{CZ10}, \cite{KW15}, \cite{KW24}, \cite{MW17}). 
Rotational symmetry for K\"ahler-Ricci solitons has also been studied (see \cite{FIK03}).

Therefore, we consider the rotational symmetry problem for gradient Yamabe solitons, which is analogous to Perelman's conjecture for Ricci solitons.
The problem can be stated as follows:

{\bf 
Under what conditions do nontrivial complete gradient shrinking or steady Yamabe solitons with positive (scalar) curvature or with bounded (scalar) curvature become rotationally symmetric? 
}

For these problems, we can provide some examples that are not rotationally symmetric.
\begin{example}\label{eg1}
Let $(N^{n-1},~\bar g)$ be an $(n-1)$-dimensional complete Riemannian manifold with constant positive scalar curvature $\bar R$.
Then for any $\alpha \in \mathbb{R}$, $(M,g,F,\lambda>0)=(\mathbb{R}\times N^{n-1}, dr^2+\frac{\bar R}{\lambda}\bar g, \sqrt{\frac{\bar R}{\lambda}} r+\alpha,\lambda)$ is an $n$-dimensional complete shrinking gradient Yamabe soliton with $R=\lambda.$
Therefore, for higher dimensions, it is possible to construct examples. 
\end{example}

Therefore, for shrinking solitons, if rotational symmetry can be shown under the condition $R \geq \lambda$ on $M$ and $R > \lambda$ at some point, this assumption is optimal.

So far, much effort has been devoted to these problems.
Since the equations for Yamabe solitons are superficially similar to those for Ricci solitons, the approaches that have driven progress on these problems have been inspired by work on Ricci solitons.
In 2012 and 2013, there was significant progress in this field. 
Inspired by the work of Cao and Chen on Ricci solitons \cite{CC12,CC13}, Daskalopoulos and Sesum \cite{DS13} provided answers to these problems under the assumptions of (i) local conformal flatness and (ii) positive sectional curvature. 
Subsequently, Cao, Sun, and Zhang \cite{CSZ12} achieved a significant breakthrough in the study of nontrivial complete gradient Yamabe solitons. They provided a structure theorem for gradient Yamabe solitons and answers to these problems under the assumptions of (i) local conformal flatness and (ii) positive scalar curvature. 
Catino, Mantegazza, and Mazzieri \cite{CMM12} considered an extension of Yamabe solitons and also achieved a significant breakthrough by removing the assumption of local conformal flatness, previously believed impossible to eliminate. 
They also relaxed the requirement of positive sectional curvature by assuming instead that the Ricci tensor is nonnegative everywhere and positive at some point, and applied their results to $k$-Yamabe solitons. 
For the steady case, the problem was solved by C. He \cite{He11}. (Note that in \cite{DS13} it was shown that certain steady and shrinking gradient Yamabe solitons satisfy $R \geq 0$. However, C. He proceeds with his arguments assuming that $R \geq 0$ in the general steady case. Therefore, his main theorem may require at least the assumption that $R \geq 0$.)
The remaining problems concern shrinking and expanding gradient Yamabe solitons. In 2013, Chu and Wang attempted to resolve the problem for the shrinking case (Theorem 1.5 in \cite{CW13}). However, the proof of Theorem 1.5, Case (1) contains a critical gap. 
Through personal communication with Yawei Chu, it was mutually confirmed that the above problem remains a challenging open question.

Despite these advances, no significant progress has been made on these problems for over ten years. The reason for this is as follows. 
As discussed in \cite{CMM12} and as we will discuss later, the core issue lies in the assumption of bounded scalar curvature from below. Furthermore, for Ricci solitons, since reducing the problem to local conformal flatness did not fundamentally resolve Perelman's conjecture, a similar reduction is unlikely to succeed here. 
In fact, since Yamabe solitons provide less information than Ricci solitons, it is virtually impossible that a reduction to the vanishing Weyl tensor, which was unsuccessful even for Ricci solitons, would succeed here. 
Additionally, in \cite{Maeta21}, the author established rotational symmetry for conformal gradient solitons ($\nabla\nabla F = \varphi g$), a generalization of gradient Yamabe solitons (which correspond to the case $\varphi = R - \lambda$), under the assumptions of local conformal flatness and positive scalar curvature. This implies that local conformal flatness is too restrictive an assumption when considering the rotational symmetry of gradient Yamabe solitons.

To overcome this difficulty, we avoid treating Yamabe solitons simply as analogous to Ricci solitons.
 Following the structure theorem of Cao, Sun, and Zhang \cite{CSZ12}, complete gradient Yamabe solitons can be classified into two types.
One of these is rotationally symmetric, while the other is of a warped product form $\mathbb{R}\times_{|\nabla F|}N$, where $N$ is a general Riemannian manifold and $\nabla F$ is the gradient of $F$.
However, in general, it is clearly difficult to analyze the warped product $\mathbb{R}\times_{|\nabla F|}N$ (see, for example, \cite{Petersen16}). 
Specifically, we need to show that $N$ is a sphere. 
However, the soliton equation and the assumptions of these problems provide us with information only about the scalar curvature.
To overcome this difficulty, we reduce the problem to the nonexistence of solutions to a certain ODE, using the structure of the warped product and the soliton equation.
Although finding solutions to this problem directly is challenging, through a very detailed analysis, we establish a nonexistence theorem.

Finally, we completely resolve these problems without imposing additional assumptions.

\begin{theorem}\label{main1}
Let $(M,g,F,\lambda)$ be a nontrivial, complete, shrinking gradient Yamabe soliton with $R\geq\lambda$. 
Then, $(M,g,F,\lambda)$ is either
\begin{enumerate}
\item[$(1)$]
 rotationally symmetric and $([0,+\infty),dr^2)\times_{{F'(r)}^2}(\mathbb{S}^{n-1},{\bar g}_{S})$, or 
\item[$(2)$] a product Riemannian manifold 
$
(\mathbb{R},dr^2)\times_{\frac{\bar R}{\lambda}} \left(N^{n-1},\bar g\right),
$
where the scalar curvature $\bar R$ of $N$ is positive constant and $F'(r)=\sqrt{\frac{\bar R}{\lambda}}$.
\end{enumerate}
\end{theorem}

Therefore, we conclude that any nontrivial, complete, shrinking gradient Yamabe soliton with $R\geq\lambda$ on $M$ and $R>\lambda$ at some point is rotationally symmetric.
As mentioned above, the assumption is optimal. 

By a similar argument, one can show that any nontrivial, complete, steady gradient Yamabe soliton with $R\geq0$ on $M$ and $R>0$ at some point is rotationally symmetric (Remark \ref{main2}), which was shown by C. He \cite{He11}.

\section{Preliminaries and the proof}\label{sectionP}

An $n$-dimensional Riemannian manifold $(M^n,g)$ is called a gradient Yamabe soliton if there exists a smooth function $F$ on $M$ and a constant $\lambda\in \mathbb{R}$ such that 
$\nabla \nabla F=(R-\lambda)g.$
We denote a gradient Yamabe soliton as $(M^n, g, F, \lambda)$.
If $F$ is constant, then $M$ is called trivial.
If $\lambda>0$, $\lambda=0$, or $\lambda<0$, then the Yamabe soliton is called shrinking, steady, or expanding, respectively.
The Riemannian curvature tensor is defined by 
\[
R(X,Y)Z=-\n_X\n_YZ+\n_Y\n_XZ+\n_{[X,Y]}Z,
\]
for $X,Y,Z\in \mathfrak{X}(M)$, where $\nabla$ is the Levi-Civita connection of $M$. The Ricci tensor $R_{ij}$ is defined by 
$R_{ij}=R_{ipjp},$ where $R_{ijk\ell}=g(R({\partial_i,\partial_j})\partial_k,\partial_\ell).$

As mentioned above, Cao, Sun, and Zhang showed the structure theorem for Yamabe solitons (see also \cite{Tashiro65} and \cite{CMM12}).

\begin{theorem}[\cite{Tashiro65}, \cite{CSZ12}, \cite{CMM12}]\label{CSZ}
Let $(M,g,F)$ be a nontrivial, complete, gradient Yamabe soliton. Then $(M,g)$ is a complete warped product manifold and must take one of the two forms:

$(1)$ 
$$([0,\infty),dr^2)\times_{{F'(r)}^2}(\mathbb{S}^{n-1},{\bar g}_{S}),$$
where $\bar g_{S}$ is the round metric on $\mathbb{S}^{n-1},$ or

$(2)$ 
$$(\mathbb{R},dr^2)\times_{{F'(r)}^2} \left(N^{n-1},\bar g\right),$$
where the scalar curvature $\bar R$ of $N$ is constant. If the scalar curvature of $M$ is nonnegative, then $\bar R>0$ or $R=\bar R=0.$ 
\end{theorem}

\begin{remark}\label{depends r}
~
\begin{enumerate}
\item[$(1)$]
For a conformal soliton, that is, a Riemannian manifold $(M,g,F,\varphi)$ that satisfies the condition $\nabla\nabla F=\varphi g$, where $\varphi\in C^\infty(M)$, Tashiro \cite{Tashiro65} provided the structure theorem, and Catino, Mantegazza, and Mazzieri also provided the structure theorem in \cite{CMM12} $($see also \cite{Maeta21}$).$ 
Such manifolds were also studied by Cheeger and Colding \cite{CC96}.

\item[$(2)$] From the proof of Theorem 1.1 in \cite{Maeta21}, one obtains that if $\Ric(\nabla F,\nabla F)\geq0$ holds everywhere and is strictly positive at some point, then such a complete Riemannian manifold $(M,g,F, \varphi)$ with a nonconstant function $F$ is rotationally symmetric. 
Therefore, in the problem of the Yamabe soliton version of Perelman's conjecture, working under scalar curvature assumptions is of fundamental importance.

\item[$(3)$] It is shown that $F$ depends only on $r$, and in Case (2) of Theorem \ref{CSZ}, without loss of generality, we can assume that $\rho(r) = F'(r) > 0$ on $\mathbb{R}$ (cf. \cite{CSZ12} see also \cite{Maeta21}).
\end{enumerate}
\end{remark}

The following lemma will play a fundamental role in the proof of our main result later.

\begin{lemma}\label{zero}
Let $a_i$ $(i=1,2,3,4)$ be constants. Suppose that $\rho: \mathbb{R} \to \mathbb{R}$ is a smooth positive function satisfying the differential equation
\begin{equation} \label{eq:rho_ode}
\rho\rho'' = a_1- a_2\rho^2 - a_3^2\rho'\rho^2 - a_4^2\rho'^2.
\end{equation}
If $\rho$ satisfies either $\rho' \geq 0$ or $\rho' \leq 0$ on $\mathbb{R}$, and there exists a point $x_0 \in \mathbb{R}$ such that $\rho'(x_0) = 0$, then $\rho' \equiv 0$ on $\mathbb{R}$.
\end{lemma}

\begin{proof}
Since $\rho'$ does not change sign on $\mathbb{R}$ and $\rho'(x_0) = 0$, the function $\rho'$ attains its global minimum or maximum at $x_0$, which yields $\rho''(x_0)=0$.
Evaluating the given differential equation \eqref{eq:rho_ode} at $x = x_0$, and using $\rho'(x_0) = 0$ and $\rho''(x_0) = 0$, we obtain
\[
0 = a_1 - a_2\rho(x_0)^2.
\]
Let $\rho_0 = \rho(x_0) > 0$. Then we have $a_1 - a_2\rho_0^2 = 0$. We can rewrite \eqref{eq:rho_ode} in the standard form for a second order ordinary differential equation:
\[
    \rho'' = \rho^{-1} \left( a_1- a_2\rho^2 - a_3^2\rho'\rho^2 - a_4^2\rho'^2 \right).
\]
Since $\rho$ is strictly positive on $\mathbb{R}$, the right-hand side is smooth with respect to $\rho$ and $\rho'$, which satisfies the local Lipschitz condition. 
Consider the constant function $\tilde{\rho}(x) \equiv \rho_0$. It is clear that $\tilde{\rho}' \equiv 0$ and $\tilde{\rho}'' \equiv 0$. Since $a_1 - a_2\rho_0^2 = 0$, the constant function $\tilde{\rho}(x)$ is a valid solution to the differential equation.

Both $\rho(x)$ and $\tilde{\rho}(x)$ satisfy the same initial conditions at $x_0$:
\[
\rho(x_0) = \tilde{\rho}(x_0) = \rho_0, \quad \rho'(x_0) = \tilde{\rho}'(x_0) = 0.
\]
By the uniqueness theorem for solutions of ordinary differential equations, they must coincide on their maximal interval of existence. Therefore, $\rho(x) = \tilde{\rho}(x) \equiv \rho_0$ for all $x \in \mathbb{R}$, which implies $\rho' \equiv 0$ on $\mathbb{R}$.
\end{proof}

\begin{lemma}\label{equivlambda}
Let $(M,g,F,\lambda)$ be a nontrivial, complete, gradient Yamabe soliton such that 
$M=\mathbb{R}\times N^{n-1}$ and $g=dr^2+{{F'(r)}^2} \bar g$.
If $R\geq\lambda$ (resp. $R\leq\lambda$) on $M$, then $R>\lambda$ or $R\equiv\lambda$ (resp. $R<\lambda$ or $R\equiv\lambda$) on $M$.
\end{lemma}

\begin{proof}
In this case, for $a,b,c,d=2,3,\cdots,n,$ the curvature tensors are obtained as follows:
\begin{align}\label{RT1}
R_{1a1b}&=-F'F'''{\bar g}_{ab},\quad R_{1abc}=0,\\
R_{abcd}&=(F')^2{\bar R}_{abcd}+(F'F'')^2(\bar g_{ad}\bar g_{bc}-\bar g_{ac}\bar g_{bd}),\notag
\end{align}
\begin{align*}
R_{11}=&-(n-1)\frac{F'''}{F'},\quad 
R_{1a}=0,\\
R_{ab}=&\bar R_{ab}-((n-2)(F'')^2+F'F''')\bar g_{ab},\notag
\end{align*}
\begin{align}\label{RT3}
R=(F')^{-2}\bar R-(n-1)(n-2)\Big(\frac{F''}{F'}\Big)^2-2(n-1)\frac{F'''}{F'},
\end{align}
where the curvature tensors with bar are the curvature tensors of $(N,\bar g)$.
By \eqref{RT3} and the soliton equation $R-\lambda=\rho'$, one has the following ODE.
\[
\lambda\rho^2+\rho'\rho^2+(n-1)(n-2)\rho'^2+2(n-1)\rho\rho''=\bar R,
\]
where $\rho=F'$. 
By Lemma \ref{zero}, we complete the proof.
\end{proof}

\begin{proof}[Proof of Theorem \ref{main1}]

To show the rotational symmetry of $M$, we only need to consider $(2)$ of Theorem \ref{CSZ}. 
By Lemma \ref{equivlambda}, we consider $R>\lambda$.
The problem can be reduced to the nonexistence of solutions for a certain ODE. By the same argument as in Lemma \ref{equivlambda}, we have 
\begin{equation}\label{key2}
\lambda\rho^2+\rho'\rho^2+(n-1)(n-2)\rho'^2+2(n-1)\rho\rho''=\bar R,
\end{equation}
where $\rho=F'$. 
By differentiating both sides of \eqref{key2}, we have
\begin{equation}\label{key3}
2\lambda\rho\rho'+\rho''\rho^2+2\rho\rho'^2+2(n-1)^2\rho'\rho''+2(n-1)\rho\rho'''=0.
\end{equation}
Since $\rho'>0$, $\rho$ is monotonically increasing.

Furthermore, one can show that $\rho$ tends to infinity. 
\begin{claim}\label{claimshr1}
$\rho(r)\nearrow+\infty$ as $r\nearrow+\infty$.
\end{claim}
\begin{proof}[Proof of Claim \ref{claimshr1}]
Assume that $\rho$ is bounded from above, that is, $F'=\rho\nearrow c$ for some positive constant $c$ as $r\nearrow+\infty$. 

If $\rho''\geq0$ on $\mathbb{R}$, since $\rho$ is bounded from above, we have that $\rho$ is constant, which is a contradiction, because $\rho'>0.$

If $\rho''\leq0$ on $\mathbb{R}$, then $\rho$ is constant because $\rho>0$. However, this cannot occur, because $\rho'>0$. 

Hence, there exist some intervals $\Omega_-$ on which $\rho''<0$ and $\Omega_+$ on which $\rho''>0$. 
There exists $r_1\in\mathbb{R}$ such that $\rho''(r_1)=0.$
By \eqref{key3}, we have $\rho'''(r_1)<0$. The same argument shows that $\rho'$ is weakly decreasing after the point $r_1$, that is, on $(r_1,+\infty)$. Without loss of generality, we may assume that $r_1$ is the first point such that $\rho''=0$, and therefore $\Omega_+=(-\infty,r_1)$. The maximum point of $\rho'$ is $r_1$. If $\rho'''<0$ on $\Omega_+$, then $\rho'$ is strictly concave on a semi-infinite interval, which contradicts $\rho'>0$ everywhere. Therefore, $\rho'''\geq0$ at some point $r_0\in\Omega_+$. However, by \eqref{key3}, we have a contradiction.

Therefore, $\rho(r)$ tends to infinity as $r\nearrow+\infty.$
\end{proof}

The equation \eqref{key2} is an autonomous second-order equation and can be reduced to a first-order equation by using $\rho$ as a new independent variable.
If $\rho'=G(\rho)$, then $\rho''=\dot G(\rho) G,$ and one has
\begin{equation}\label{1st}
\lambda\rho^2+G\rho^2+(n-1)(n-2)G^2+2(n-1)\rho\dot GG=\bar R.
\end{equation}
By differentiating the equation, one has
\begin{equation}\label{2nd}
2\lambda\rho+\dot G\rho^2+2\rho G+2(n-1)^2\dot G G+2(n-1)\rho\ddot G G+2(n-1)\rho \dot G^2=0.
\end{equation}
Assume that $\dot G>0$ at some point $\rho_0\in (0,+\infty)$, that is, $\dot G>0$ on some open interval $\Omega=(\rho_1,\rho_2)(\ni\rho_0)$. If $\Omega=(\rho_1,+\infty)$, then by \eqref{1st}, 
$\lambda\rho^2+G\rho^2+(n-1)(n-2)G^2<\bar R$ on $(\rho_1,+\infty).$
However, the left-hand side tends to infinity as $\rho\nearrow +\infty$ by Claim \ref{claimshr1}, which cannot occur.
Thus, one can assume that $\dot G=0$ at $\rho_2$.
Then by $\eqref{2nd}$, $2\lambda\rho_2+2\rho_2 G+2(n-1)\rho_2\ddot GG=0$ at $\rho_2$. 
Hence, $\ddot G<0$ at $\rho_2$, and $G$ is monotonically decreasing on $(\rho_2,\rho_3)$ for some $\rho_3$. By iterating the same argument, one can extend $\rho_3$ to $+\infty$, that is, $G$ is monotonically decreasing on $(\rho_2,+\infty)$. 
Hence, if such an open interval $\Omega$ exists, it must be $(0,\rho_2)$ and $\rho_2$ is the maximum point of $G$. 
Let $p\in\mathbb{R}$ be the maximum point of $\rho'$. Since $\rho'>0$, there exists a point $q\in (-\infty, p)$ such that $\rho''(q)>0$ and $\rho'''(q)=0$. However, by \eqref{key3}, one has
$2\lambda\rho\rho'+\rho''\rho^2+2\rho \rho'^2+2(n-1)^2\rho'\rho''=0$ 
at $q$, which cannot occur. We finally obtain $\Omega=\emptyset.$ 

Therefore, one has $0\geq\dot G=\frac{\rho''(r)}{\rho'(r)}$ for every $r\in\mathbb{R}$, and $\rho''\leq 0$ on $\mathbb{R}$. 
Since the positive smooth function $\rho$ is concave on $\mathbb{R}$, $\rho$ is constant, which cannot occur. 
\end{proof}

\begin{remark}\label{main2}
By an argument similar to the proof of Theorem \ref{main1}, we can also establish the following rigidity result for steady solitons: {\it Any nontrivial, complete, steady gradient Yamabe soliton whose scalar curvature is nonnegative everywhere and is strictly positive at some point is rotationally symmetric.} This approach provides an alternative proof of the result originally obtained by He \cite{He11}.
\end{remark}

\noindent
{\bf Acknowledgements.}~\\
The author would like to express his gratitude to Giovanni Catino and Takumi Yokota for their useful comments, and to Ken Shirakawa and Yawei Chu for valuable discussions. 
~\\




\bibliographystyle{amsbook}

\end{document}